\documentclass[11pt]{amsart}
\usepackage{amsmath,amsfonts,amssymb,amsxtra,latexsym,amscd,enumerate,amsthm}

\usepackage{latexsym}
\usepackage{epsfig}
\usepackage{verbatim}

\def\g{\gamma}
\def\a{\alpha}

\def\d{\delta}

\def\p{\Phi}
\def\f{\varphi}
\def\ep{\epsilon}

\def\ll{\lambda}

\def\o{\omega}

\def\cN{\mathcal{N}}
\def\R{\mathbb{R}}

\def\Z{\mathbb{Z}}

\def\D{\mathcal{D}}
\def\P{\mathbb{P}}
\def\p{\mathcal{P}}
\def\l{\mathcal{L}}
\def\I{\mathcal{I}}
\def\J{\mathcal{J}}
\def\LL{\mathcal{L}}

\def\T{\mathcal{T}}
\def\F{\mathcal{F}}
\def\A{\mathcal{A}}

\def\TT{\mathbb{T}}
\def\N{\mathbb{N}}
\def\W{\mathcal{W}}
\def\Z{\mathbb{Z}}
\def\beq{\begin{equation}}
\def\eeq{\end{equation}}

\def\beq{\begin{equation}}
\def\eeq{\end{equation}}

\newtheorem{d0}{Definition}

\newtheorem{l1}{Lemma}
\newtheorem{p1}{Proposition}

\begin{document}
\title[]{On the Boundedness of the Carleson Operator near $L^1$}

\author{Victor Lie}

\date{\today}

\address{Department of Mathematics, Princeton, NJ}

\email{vlie@math.princeton.edu}
\address{Institute of Mathematics of the Romanian Academy, Bucharest, RO
70700 \newline \indent  P.O. Box 1-764}

\keywords{Time-frequency analysis, Carleson's Theorem.}

\maketitle

\begin{abstract}
Based on the tile discretization elaborated in \cite{lv3}, we develop a Calderon-Zygmund type decomposition of the Carleson
operator. As a consequence, through a unitary method that makes no use of extrapolation techniques,
we recover the previously known results regarding the largest rearrangement
invariant space of functions with almost everywhere convergent Fourier series.
\end{abstract}

\section{Introduction}

In this paper we analyze some aspects concerning the behavior of the Carleson operator near $L^1$. Elaborating on an idea\footnote{This approach has as a consequence the removal of the exceptional sets in the tile decomposition and thus gives direct strong $L^2$ bounds for the Carleson operator. Also this is one of the key ingredients in providing the full range for the $L^p$ bounds ($1<p<\infty$) of the Polynomial Carleson operator.}
introduced in \cite{lv3}, we design a Calderon-Zygmund type decomposition of the Carleson operator which, besides the direct interest in our problem, may prove useful in other related topics. In particular, through this technique we are able to encompass the previously known results regarding the problem of the largest rearrangement invariant space of (integrable) functions for which the \textit{a.e.} convergence of the Fourier Series holds. The relevant point here though is not the ability of reproving these results but rather the existence of a method that avoids the limitations of the extrapolation techniques - the main ingredient on which all the previous results rely on.

As we will see, most of the difficulty and interest resides in the tile decomposition of the Carleson operator. Once this decomposition is achieved, everything else follows naturally. In the present paper we are focussing on the method rather than on obtaining the best possible space
on which the Carleson operator is finitely almost everywhere. With respect to the latter the best current result belongs to Arias de Reyna (\cite{Ar}).

A significant improvement of his result seems to require an original idea, this paper setting just the first step in opening up a new direction\footnote{See also 3) in the Remarks section.} of investigation.

These being said, let us state the precise theme of interest regarding the behavior of the Fourier Series near $L^1$:
$\newline$

\textbf{Open problem.} \textit{What is the largest Banach rearrangement invariant space $(Y,\,\|\cdot\|_{Y})$ with $Y\subseteq L^1(\TT)$ for which
 the Carleson operator defined\footnote{We set $D'(\TT)$ the class of distributions supported on the torus.} by $$T:\,C^{\infty}(\TT)\,\mapsto\,\D'(\TT)$$ with
\beq\label{carl}
Tf(x):=\sup_{N\in\N}\left|\int_{\TT}e^{i\,N\,(x-y)}\,\cot(x-y)\,f(y)\,dy\right|\:,
\eeq
obeys the relation
\beq\label{carlb}
\|Tf\|_{1,\infty}\lesssim\|f\|_{Y}\:\:\:\:\:\forall\:\:f\in Y\:?
\eeq}

\textbf{Observation.} \textit{Notice that from Stein's maximal principle (\cite{s1}) the above open problem is equivalent to asking for the largest Banach rearrangement invariant space $(Y,\,\|\cdot\|_{Y})$ with $Y\subseteq L^1(\TT)$ for which the partial Fourier Series $\{S_n f\}_n$ have the property}
 \beq\label{pointconv}
S_n f(x)\,\stackrel{n \rightarrow \infty}{\longrightarrow}\,f(x)\:\:a.e.\:\:x\in\TT\:\:\:\:\:\forall\:\:f\in Y.
\eeq

 Significant literature has been written on this subject. The major breakthrough was made by Carleson (\cite{c1}), who showed that $L^2(\TT)\subset Y$. Later, Hunt (\cite{hu}) extended this result by showing that $L^p(\TT)\subset Y$ for any $1<p<\infty$. A new influential proof of Carleson's result\footnote{More recently (\cite{lt3}), Lacey and Thiele, combining ideas from both \cite{c1} and \cite{f}, provided a third approach to Carleson's theorem on the pointwise convergence of the Fourier Series.} was given by Fefferman in \cite{f}. From this on, the problem evolved at a slower rate towards the limiting index $p=1$:
 Sj\"olin, proved in \cite{sj3}, that one may take in \eqref{carlb} $Y=L (\log L)^2$ (if one requires strong $L^1$ bounds in \eqref{carlb}) or even $L\,\log L\,\log\log L$ for just $L^{1,\infty}$ bounds. Next, Soria (\cite{So1},\cite{So2}) constructed a larger space\footnote{See the Appendix for definition.} $B_{\f_1}^*\subset Y$. In \cite{An}, Antonov showed that \eqref{carlb} holds for $Y=L\log L \log\log\log L$. Finally, combining elements from Antonov's and Soria's approaches with techniques on logconvex quasi-Banach spaces, Arias de Reyna (\cite{Ar}) proved that $QA\subset Y$, where $QA$ is a quasi-Banach space described in the Appendix.

 The following chain of inclusions\footnote{Notice that there is no order relation between the spaces $B_{\f_1}^*$ and $L\,\log L\,\log\log\log L$ (for more details see \cite{Ar}).} holds:
\beq\label{relspaces}
L (\log L)^2\subsetneq L\,\log L\,\log\log L\subsetneq B_{\f_1}^*,\,L\,\log L\,\log\log\log L \subsetneq QA\subsetneq L \log L\:.
\eeq

The main result in this paper is given by
$\newline$

\noindent\textbf{Theorem.} \textit{There exists a partition of the family of tiles
$$\P=\bigcup_{n\in\N}\P_n\,,$$
with $$T=\sum_{n} T^{\P_n}$$ such that for each $n\in\N$ we have}

a) \textit{Given $f\in L^1$ nonzero, there exists a further decomposition\footnote{This second decomposition depends on $f$.}
$$\P_n=\bigcup_{\a\in\Z}\P_n^{\a}\,,$$
such that for any $n\in\N$ we have\footnote{Here, if $J=(c-\frac{|J|}{2}, c+\frac{|J|}{|2|})$ is any given interval, we use the standard notation
$b\,J$ ($b>0$) to designate the interval $(c-\frac{b\,|J|}{2}, c+\frac{b\,|J|}{|2|})$. Moreover, if $I=\bigcup_{n\in\N} J_n$ with each $J_n$ an interval, then we set $b\,I:=\bigcup_{n\in\N} b\,J_n$.}
\beq\label{suport}
\textrm{supp}\,T^{\P_n^{\a}}\subseteq 100\,\{Mf>2^{-\a}\}\,,
\eeq
\beq\label{main}
\|T^{\P_n^{\a}}\,f\|_1\lesssim 2^{-\a}\,|\{Mf>2^{-\a}\}|\,,
\eeq
where here $M$ stands for the dyadic Hardy-Littlewood maximal function.}

\textit{In particular, we deduce that}
\beq\label{llogl}
\|T^{\P_n}\,f\|_{1}\lesssim \|f\|_{L\,\log L}\:.
\eeq

b) \textit{For $f\in L^1$ the following holds:}
\beq\label{l1}
\|T^{\P_n}\,f\|_{1,\infty}\lesssim \|f\|_{1}\:.
\eeq

c) \textit{If  $1<p<\infty$ and} $p^{*}=\textrm{min}\{p,p'\}$ \textit{with $p'$ the H\"older conjugate of $p$, then there exists
an absolute constant $\d>0$ such that}
\beq\label{lp}
\|T^{\P_n}\,f\|_{p}\lesssim_p 2^{-\d\,n/p^*}\,\|f\|_{p}\:.
\eeq

d) \textit{If $f\in L^p$ with $1<p\leq\infty$ then}
\beq\label{linterp}
\|Tf\|_{1,\infty}\lesssim_{p} \|f\|_1\,\log\frac{e\,\|f\|_{p}}{\|f\|_1}\:.
\eeq
$\newline$

As a \textsf{direct} application of our Theorem we have the following
$\newline$

\noindent\textbf{Corollary.} \textit{The following are true:}

1) (\textsf{Carleson-Hunt}, \cite{c1}, \cite{hu}) \space     \textit{$\|Tf\|_p\lesssim_p\|f\|_p$ for any $1<p<\infty$.}

2) (\textsf{Sj\"olin}, \cite{sj3}) \space  $\|Tf\|_1\lesssim\|f\|_{L(\log L)^2}\;.$

3) \textit{For $E\subseteq [0,1]$ measurable $\|T\chi_{E}\|_{1,\infty}\lesssim |E|\,\log\frac{e}{|E|}\;.$}

4) (\textsf{Sj\"olin}, \cite{sj3}) \space  $\|Tf\|_{1,\infty}\lesssim\|f\|_{L\log L\log\log L}\;.$

5) (\textsf{Arias de Reyna}, \cite{Ar}) \space  $\|Tf\|_{1,\infty}\lesssim \|f\|_{QA}\;.$

\textit{In particular 5) also implies}

6) (\textsf{Soria}, F., \cite{So1},\cite{So2}) \space  $\|Tf\|_{1,\infty}\lesssim \|f\|_{B^{*}_{\f_1}}\;.$

7) (\textsf{Antonov}, \cite{An}) \space  $\|Tf\|_{1,\infty}\lesssim \|f\|_{L\log L\log\log\log L}\;.$
$\newline$

\noindent\textbf{Comment.} \textit{In fact, as a consequence of d) in our Theorem, we obtain that for $1<p\leq\infty$ one has}
$$\|Tf\|_{1,\infty}\lesssim_p \|f\|_{QA_p}\;,$$
\textit{where $QA_p$ is the quasi-Banach space defined by
$$QA_p:=\{f:\:\TT\mapsto C\,|\,f\:\textrm{measurable},\:\|f\|_{QA_p}<\infty \}\,$$
with
\beq\label{W}
\|f\|_{QA_p}:=\inf\left\{\sum_{j=1}^{\infty}(1+\log j)\|f_j\|_1\,\log\frac{e\,\|f_j\|_{p}}{\|f_j\|_1}\:\:\left|\right.\:\:
\begin{array}{cl}
f=\sum_{j=1}^{\infty}f_j,\:\\
\sum_{j=1}^{\infty}|f_j|<\infty\:\textrm{a.e.}
\end{array}
\right\}\;.
\eeq}
\textit{But as it turns out, based on an observation of Louis Rodriguez-Piazza kindly provided to me by Arias de Reyna, the spaces $QA_p$
are equivalent in the sense that $$\|f\|_{QA_p}\approx_{p}\|f\|_{QA_\infty}\:.$$
Notice that $QA_{\infty}$ is the original space $QA$ defined in \cite{Ar}; for an interesting study concerning the properties of the space $QA$ one should consult \cite{CMR}.}
$\newline$

\noindent\textbf{Acknowledgments:} We thank Arias de Reyna for reading the manuscript and supplying with useful comments.

\section{Discretization of the operator}

In this section we decompose the operator $T$ in components $\{T_P\}_P$ which are ``well" time-frequency localized.
We follow the procedure in \cite{f}.

Let $T$ be the Carleson operator defined by
$$Tf(x):=\sup_{N\in\R}\,|\int_{\TT} \frac{1}{x-y}\,e^{i\,N\,(x-y)}\,f(y)\,dy| \:,$$
which after linearization becomes
$$Tf(x)=\int_{\TT} \frac{1}{x-y}\,e^{i\,N(x)\,(x-y)}\,f(y)\,dy\:,$$
where $N$ is some arbitrary measurable function (which from now will be fixed).

Choose now $\psi$ an odd
$C^{\infty}$ function such that
$$\operatorname{supp}\:\psi\subseteq\left\{y\in \R\:|\:2<|y|<8\right\}$$ which has the property
$$\frac{1}{y}=\sum_{k\geq 0} \psi_k(y)\:\:\:\:\:\:\:\:\:\forall\:\:0<|y|<1\:,$$
where by definition $\psi_k(y):=2^{k}\psi(2^{k}y)$ (with $k\in \N$).

Thus, we have that
$$Tf(x)=\sum_{k\geq 0}T_{k}f(x):=\sum_{k\geq 0}\int_{\TT}e^{i\,N(x)\,y}\psi_{k}(y)f(x-y)dy\:.$$

Take the canonical dyadic grid on $\TT$ (time grid) - denoted by $\D_{\T}$ and the corresponding canonical dyadic grid on $\R$
(frequency grid) - denoted with $\D_{\F}$. A tile $P$ will consist from a tuple\footnote{Here we abuse the language and refer
to $\D_{\T}$ also as the collection of all dyadic intervals for the specified grid. Same for $\D_{\F}$.} $[\o, I]\in\D_{\F}\times\D_{\T}$
with the property that $|\o|=|I|^{-1}$. The collection of all such tiles will be denoted with $\P$.
Further, for each $P=[\o,I]\in\P$ we set $E(P):=\left\{x\in I\:|\:N(x)\in \o\right\}$.

With this being said, for $|I|=2^{-k}$ ($k\geq0$) and $P=[\o,I]\in\P$ we define the operators
$ T_P$ on $L^2(\TT)$ by
$$T_{P}f(x)=\left\{\int_{\TT}e^{i\,N(x)\,y}\psi_{k}(y)f(x-y)dy\right\}\chi_{E(P)}(x)\:.$$

Notice that the Carleson operator obeys
\beq\label{discret}
Tf(x)=\sum_{P\in\P} T_P f(x)\,.
\eeq
\indent Finally, whenever $\p\subseteq\P$ is a family of tiles we set
$$T^{\p}:=\sum_{P\in\p}T_P\:.$$

\section{The proof of the Corollary.}

\subsection{The proof of 1).} We want to show that for $1<p<\infty$
\beq\label{CH}
\|Tf\|_p\lesssim_p\|f\|_p\:.
\eeq
This is a trivial application of statement c) in our Theorem. Indeed, we have
$$\|Tf\|_p\leq\sum_{n\in N} \|T^{\P_n}\,f\|_{p}\lesssim_p\sum_{n\in\N} 2^{-\d n/p^{*}}\,\|f\|_{p}\lesssim\|f\|_p\:.$$

\subsection{The proof of 2).} We want to show that
\beq\label{Sj}
\|Tf\|_1\lesssim\|f\|_{L(\log L)^2}\;.
\eeq
We will use the following decomposition
for each $l\in\Z$ define $$Q_{l}:=\{x\in\TT\,|\,|f(x)|\in [2^{l},\,2^{l+1})\}\:.$$ Then we have that
\beq\label{llogl}
\|f\|_{L (\textrm{log}L)^2}\sim \sum_{l\in\Z} 2^{l}\,|Q_l|\,(\log \frac{1}{|Q_l|})^2\:.
\eeq
Thus, using duality, for proving \eqref{Sj} will be enough to show
\beq\label{enough}
\int_{Q_l} |T^{*}g| \lesssim |Q_l|\,\large(\log \frac{1}{|Q_l|}\large)^2\, \|g\|_{\infty}\:.
\eeq
Taking $f=\chi_{Q_l}=\chi_{Q}$  in the dual statements of a) and c) in our Theorem  we deduce
$$\int_{Q} |{T^{\P_n}}^{*}g| \lesssim |Q|\,\log \frac{1}{|Q|}\, \|g\|_{\infty}\:\:\textrm{and}\:\:
\int_{Q} |{T^{\P_n}}^{*}g| \lesssim |Q|^{1/2}\,2^{-n\,\d}\, \|g\|_{\infty}\;.$$
Thus
$$\int_{Q} |{T^{\P_n}}^{*}g|\lesssim \sum_n\min\{|Q|\,\log \frac{1}{|Q|},\,|Q|^{1/2}\,2^{-n\,\d}\}\|g\|_{\infty}\lesssim
|Q|\,\large(\log \frac{1}{|Q|}\large)^2\, \|g\|_{\infty}\:.$$

\subsection{The proof of 3).} We are interested in
\beq\label{log}
\|T\chi_{E}\|_{1,\infty}\lesssim |E|\,\log\frac{e}{|E|}\;.
\eeq
Just apply d) with $f=\chi_E$.

\subsection{The proof of 4).} Of course, this result is implied by the claims 5)-7). Still, we think it is worth providing a different approach to this problem, one that isolates a relevant idea in Sj\"olin's original proof and nicely adapts it in the context of our Theorem.

Our task is to prove that
\beq\label{Sjw}
\|Tf\|_{1,\infty}\lesssim\|f\|_{L\log L\log\log L}\;.
\eeq

The definition that we take here for $L\log L\log\log L$ is as in the original paper of Sj\"olin (\cite{sj3}) given by the space of functions
$f\in L^1(\TT)$ for which we have \footnote{Notice that unlike $\|f\|_{L\log L\log\log L}$,  $\|f\|_{\l\log \l\log\log \l}$ is not a norm.}
$\|f\|_{\l\log \l\log\log \l}=\int|f|\,\log_{+}|f|\,\log_{+}\log_{+}|f|<\infty\:.$

Thus, setting as before $Q_{l}:=\{x\in\TT\,|\,|f(x)|\approx 2^{l}\}$, we can always assume $Q_l=\emptyset$ for $l\leq 2$ and consequently
$$\|f\|_{\l\log \l\log\log \l}\approx\sum_{l>2} 2^l\,l\,\log l\,|Q_{l}|\;.$$

Now fix $l>2$. For $\chi_{l}=\chi_{Q_l}$ we run\footnote{We assume here that $|Q_l|\not=0$.} the tile partition described at point a) of our Theorem. Then for each $n\in\N$ there exists a decomposition of $\P_n=\bigcup_{\a\in\Z}\P_n^{\a}$ such that
\beq\label{q1}
\textrm{supp}\,T^{\P_n^{\a}}\subseteq 100\,\{M \chi_{l}>2^{-\a}\}\:\:\:\textrm{and}\:\:\:
\|T^{\P_n^{\a}}\,\chi_{l}\|_1\lesssim 2^{-\a}\,|\{M \chi_{l}>2^{-\a}\}|\,.
\eeq
Split the set $\Z=A_l\cup B_l\cup C_l$ where
$$A_l:=\{r\in N\,|\,2^{-r}<\g\, 2^{-l}\,l^{-3}\}\:,$$
$$B_l:=\{r\in N\,|\,\g\, 2^{-l}\,l^{-3}\leq 2^{-r}<\g\, 2^{-l}\}\:,$$
$$C_l:=\{r\in N\,|\,2^{-r}\geq\g\, 2^{-l}\}\:,$$
with $\g>0$ a parameter that will be chosen later.

Set $$T^{\P^{\a}}:=\sum_{n\in\N} T^{\P_n^{\a}}\;.$$
Then, based on \eqref{q1}, we have
\beq\label{q2}
\textrm{supp}\,T^{\P^{\a}}\subseteq 100\,\{M \chi_{l}\geq2^{-\a}\}\,,
\eeq
which implies that
\beq\label{suportt}
\sum_{\a\in C_l}|\textrm{supp}\,T^{\P^{\a}}|\lesssim\sum_{\a\in C_l}|\{M \chi_{l}\geq 2^{-\a}\}|\lesssim
\sum_{\a\in C_l} 2^{\a}\,|Q_l|\lesssim \g^{-1}\,2^l\,|Q_l|\,.
\eeq
Thus the set $S_l=\bigcup_{\a\in C_l}\,\textrm{supp}\,T^{\P^{\a}}$ can be excised since we have a good control on
$$\sum_{l\in\Z}|S_l|\lesssim\g^{-1}\,\sum_{l\in\Z}2^l\,|Q_l|\lesssim \g^{-1}\,\|f\|_{\l\log \l\log\log \l}\:.$$
Next, it is useful to notice that based on \eqref{q1} we have
$$\|T^{\P_n^{\a}}\,\chi_{l}\|_1\lesssim 2^{-\a}\,|\{M \chi_{l}\geq 2^{-\a}\}|\lesssim\min\{ 2^{-\a},\,|Q_l|\},$$
while based on d) in our Theorem we infer that
$$\|T^{\P_n^{\a}}\,\chi_{l}\|_2\lesssim 2^{-n\,\d}\,|Q_l|^{\frac{1}{2}}\,.$$
From these we deduce
\beq\label{palfa}
\|T^{\P^{\a}}\,\chi_{l}\|_1\lesssim\a\,\min\{ 2^{-\a},\,|Q_l|\}\:.
\eeq
Then, from \eqref{palfa}, we have that
\beq\label{B}
\sum_{\a\in B_l}\|T^{\P^{\a}}\,\chi_{l}\|_1\lesssim (l+\log\frac{1}{\g})\,\log l\,|Q_l|\,,
\eeq
and respectively
\beq\label{A}
\sum_{\a\in A_l}\|T^{\P^{\a}}\,\chi_{l}\|_1\lesssim \sum_{\a\in A_l} \a\, 2^{-\a}\lesssim \g \,(1+\log\frac{1}{\g})\, 2^{-l}\,l^{-2}\,.
\eeq
Putting together \eqref{suportt}, \eqref{B} and \eqref{A} and choosing $\g=c\,\|f\|_{\l\log \l\log\log \l}^{\frac{2}{3}}$ (with $c>0$ some large number) we have proved the following

\begin{p1}\label{Sjolin}
Let $f\in L\log L\log\log L$ with $\|f\|_{\l\log \l\log\log \l}<1$.

Then there exists $A\subseteq [0,1]$ with $|A|\leq\|f\|_{\l\log \l\log\log \l}^{\frac{1}{3}}$ and $C>0$ an absolute
constant such that
\beq\label{sjw}
\|Tf\|_{L^1(A^{c})}\leq C\,\|f\|_{\l\log \l\log\log \l}^{\frac{1}{2}}\:.
\eeq
\end{p1}

Now, by a canonical density argument, we obtain that the sequence of the partial Fourier sums $\{S_n f(x)\}_n$ converges almost everywhere
for $f\in L\log L\log\log L$. Relation \eqref{Sjw} follows from an application of Stein's maximal principle (\cite{s1}).

\subsection{The proof of 5)-8).} We will show that
\beq\label{WW}
\|Tf\|_{1,\infty}\lesssim_p \|f\|_{QA_p}\;.
\eeq

We choose\footnote{One can use a different approach to \eqref{WW} that avoids the use of Kalton's theorem on the log convexity of $\|\cdot\|_{1,\infty}$. For this, one can follow the proof of 4).} here to use the log-convexity result due to Kalton (\cite{Ka}), as further described in \cite{Ar}:
$\newline$

\noindent\textbf{Theorem.} \textit{If $\{f_j\}_j$ is a sequence of functions in $L^{1,\infty}(\TT)$, then we have}
\beq\label{K}
\|\sum_{j}f_j\|_{1,\infty}\lesssim \sum_{j}(1+\log j)\,\|f_j\|_{1\,\infty}\:.
\eeq

Take $f\in \W$ and set $f=\sum_{j=1}^{\infty}f_j$ with $f_j$ as in \eqref{W}. Apply now \eqref{K} and point
d) in our Theorem to conclude
$$\|Tf\|_{1,\infty}\lesssim \sum_{j}(1+\log j)\,\|T f_j\|_{1\,\infty}\lesssim_p\sum_{j=1}^{\infty}(1+\log j)\|f_j\|_1\,\log\frac{e\,\|f_j\|_{p}}{\|f_j\|_1}\:.$$

\section{Finding structures in our family of tiles; main definitions}

In this section we isolate the main concepts needed for further discretizing and organizing the family of tiles $\P$.
Our presentation here is based on the definitions and notations introduced in \cite{lv3}.

\begin{d0}\label{mass}(\textsf{weighting the tiles})

Let $\A$ be a (finite) union of dyadic intervals in $[0,1]$ and $\p$ be a finite family of tiles. For $P=[\o,I]\in\p$ with $I\subseteq\A$ we define the {\bf mass} of $P$ relative to the set of tiles $\p$ and the set $\A$ as being

\beq\label{v1} A_{\p,\A}(P):=\sup_{{P'=[\o',I']\in\:\p}\atop{I\subseteq
I'\subseteq \A}}\frac{|E(P')|}{|I'|}\:\left\lceil
\Delta(10P,\:10P')\right\rceil^{N} \eeq where $N$ is a fixed large
natural number.
\end{d0}

\begin{d0}\label{ord} (\textsf{ordering the tiles})

   Let $P_j=[\o_j,I_j]\in\P$ with $j\in\left\{1,2\right\}$. We say that
 $P_1\leq P_2$       iff       $\:\:\:I_1\subseteq I_2$
and  $\o_1\supseteq\o_2$. We write $P_1< P_2$ if $P_1\leq P_2$ and $|I_1|<|I_2|$.
\end{d0}

Notice that $\leq$ defines an order relation on the set $\P$.

\begin{d0}\label{tree} (\textsf{modulated/scaled (maximal) Hilbert transform - ``tree"})

 We say that a set of tiles $\p\subset\P$ is a \textbf{tree} with \emph{top} $P_0$ if
the following conditions hold:
$\newline 1)\:\:\:\:\:\forall\:\:P\in\p\:\:\:\Rightarrow\:\:\:\:2P\leq10 P_0 $
$\newline 2)\:\:\:\:\:$if $P=[\o_1, I_P]\in\p$ and $P'=[\o_2, I_P]$ such that $2P'\leq P_{0}$ then $P'\in\p$
$\newline 3)\:\:\:\:\:$if $P_1,\:P_2\: \in\p$ and $P_1\leq P \leq P_2$ then $P\in\p\:.$
\end{d0}

\begin{d0}\label{stree} (\textsf{Carleson measure relative to a tree})

 We say that a set of tiles $\p\subset\P$ is a \textbf{sparse tree} if
 $\p$ is a tree and for any $P\in\p$ we have
 \beq\label{cmstree}
 \sum_{{P'\in\p}\atop{I_{P'}\subseteq I_P}}|I_{P'}|\leq C\,|I_P|\:,
 \eeq
 where here $C>0$ is an absolute constant.
 \end{d0}

\begin{d0}\label{infforest} (\textsf{$L^{\infty}$ control over union of trees})

Fix $n\in\N$. We say that $\p\subseteq\P_n$ is an $L^\infty$-\textbf{forest} (of $n^{th}$-generation) if

i) $\p$ is a collection of separated trees, {\it i.e.}
$$\p=\bigcup_{j\in\N}\p_j\,$$
with each $\p_j$ a tree with top $P_j=[\o_j,I_j]$ and such that
 \beq\label{wellsep}
 \forall\:\:k\not=j\:\:\&\:\:\forall\:\:P\in\p_j\:\:\:\:\:\:\:\:2P\nleq 10P_k\:.
 \eeq

ii) the counting function
\beq\label{cfunct}
 \cN_{\p}(x):=\sum_{j}\chi_{I_j}(x)
 \eeq
 obeys the estimate $\|\cN_{\p}\|_{L^{\infty}}\lesssim 2^n$.

 Further, if $\p\subseteq\P_n$ only consists of sparse separated trees then
we refer at $\p$ as a \textbf{sparse $L^\infty$-forest}.
 \end{d0}

 \begin{d0}\label{forest} (\textsf{$BMO$ control over union of trees})

  A set  $\p\subseteq\P_n$ is called a  $BMO$-\textbf{forest} (of $n^{th}$-generation) or
 simply a \textbf{forest} if

 i) $\p$ may be written as
$$\p=\bigcup_{j\in\N}\p_j\,$$
with each $\p_j$ an $L^\infty$-\textbf{forest} (of $n^{th}$-generation);

 ii) for any $P\in\p_j$ and $P'\in\p_k$ with $j,k\in\N$,  $j<k$ we either have $I_P\cap I_{P'}=\emptyset$ or
 $$|I_{P'}|\leq 2^{j-k}\,|I_{P}|\:.$$

 As before, if $\p\subseteq\P_n$ only consists of sparse $L^\infty$-forest, then
we refer it as a \textbf{sparse forest}.
\end{d0}

 Notice that if $\p\subseteq\P_n$ is a forest then, due to {\it ii)} above, the counting function
 $\cN_{\p}:=\sum_{j}\cN_{\p_j}$ obeys  the estimate
 $$\|\cN_{\p}\|_{BMO_{C}}\lesssim 2^n\:.$$

\section{Discretization of the family of tiles}

\subsection{The mass decomposition - $n$ discretization}
$\newline$

In this section we partition the set $\P$ into $\bigcup_{n\in\N}\P_n$, with each $\P_n$ a $BMO-$forest.

The procedure described below is an adaptation of the one introduced by the author in \cite{lv3} for proving the $L^p $ boundedness ($1<p<\infty$) of the Polynomial Carleson operator.

We start by constructing the family $\P_1$ according to the following algorithm:
\begin{itemize}
\item Let $\p_{1}^{0,max}$ be the collection of maximal tiles $P\in\P$ with $\frac{|E(P)|}{|I|}\geq \frac{1}{2}$ and let
$\I_{1}^{0}:=\{I\,|\,P=[\o,I]\in\p_{1}^{0,max}\}$.
\item Define the counting function $\cN_{1}^{0}:=\sum_{I\in\I_{1}^{0}}\chi_{I}$ and verify that
\beq\label{defbmoc}
\|\cN_{1}^{0}\|_{BMO_{C}}:=\sup_{J\textrm{dyadic}\atop{J\subseteq [0,1]}}\,\frac{\sum_{I\subseteq J\atop{I\in\I_{1}^{0}}}|I|}{|J|}\leq 2\:.
\eeq
and hence  $\cN_{1}^{0}\in BMO_{D}(\R)$.
\item Apply the John-Nirenberg inequality
\beq\label{JN}
|\{x\in J\,|\,|\cN_1(x)-\frac{\int_{J} \cN_1^{0}}{|J|}|>\gamma\}|\lesssim |J|\,e^{-c\:\frac{\g}{\|\cN_1^{0}\|_{BMO_{D}(\R)}}}\:.
\eeq
for\footnote{Notice that $\|\cN_1^{0}\|_{BMO_{D}(\R)}\leq 2\|\cN_{1}^{0}\|_{BMO_{C}}$.} $\g> c\,\|\cN_1^{0}\|_{BMO_{C}}$ (here $c>0$ is an appropriately chosen large absolute constant) and deduce that
\beq\label{JNP}
|\{x\in J\,|\,\sum_{I\subseteq J\atop{I\in\I_{1}^{0}}}\chi_{I}(x)>\gamma\}|\lesssim |J|\,e^{-c}\:.
\eeq
\item Based on \eqref{JNP}, conclude that the set
$$A_1^{1}:=\{x\in [0,1]\,|\, \sum_{I\subseteq [0,1]\atop{I\in\I_{1}^{0}}}\chi_{I}(x)> c\,\|\cN_{1}^{0}\|_{BMO_{C}}\}$$
obeys the relation $|A_1^{1}|\leq e^{-c}$.
\item Remove from $\P$ all the tiles which have the time interval \emph{not} included in the set $A_1^{1}$. Run again the above algorithm
for the new collection $\P$. This process ends in a finite number of steps since wlog we may assume that the initial family $\P$ is finite.
\end{itemize}

This way, after the $k^{th}$ repetition of our algorithm, we have constructed the sets $A_1^{k}$, $\p_{1}^{k,max}$, $\I_{1}^{k}$ and the counting function $\cN_{1}^{k}$.

We now define the $1-$maximal set of tiles $\p_{1}^{max}:=\bigcup_{k}\p_{1}^{k,max}$, the collection of the time-intervals $\I_{1}:=\bigcup_{k}\I_{1}^{k}$ and finally the counting function $\cN_{1}:=\sum_{I\in\I_{1}}\chi_{I}$.

Notice that from the above construction we have that
\begin{itemize}
\item $\|\cN_{1}\|_{BMO_{C}}\lesssim \max_{k}\|\cN_{1}^{k}\|_{BMO_{C}}\:;$
\item for any $l<k$ we have $A_1^{k}\subset A_1^{l}$ and $|A_1^{k}|\leq e^{-(k-l)\:c}\,|A_1^{l}|\:.$
\end{itemize}

Next, define
$$\p_{1}^{0}:=\{P=[\o,I]\in\P\,|\,I\nsubseteq A_{1}^{1}\:\:\&\:\:A_{\P,[0,1]}(P)\in [2^{-1}, 2^{0})\}$$
and further, by induction, construct
$$\p_{1}^{k}:=\left\{P=[\o,I]\in\P\,|\,\begin{array}{rl} I\nsubseteq A_{1}^{k+1}\:,I\subseteq A_{1}^{k}\:\:\:\:\\
A_{\P,A_{1}^{k}}(P)\in [2^{-1}, 2^{0}) \end{array}\right\}\:.$$

Finally, setting $$\P_1:=\bigcup_{k}\p_{1}^{k}$$
we end the construction of the family of tiles having the mass of order $1$.

Now suppose that we have constructed the sets $\{\P_k\}_{k<n}$. Here is how we define the set $\P_n$.

Firstly, select the family $\p_{n}^{0,max}$ of the maximal tiles $P\in\P\setminus\bigcup_{k<n}\P_k$ with $\frac{|E(P)|}{|I_P|}\geq 2^{-n}$.
Then, collect the time-intervals of these maximal tiles into the set $\I_{n}^{0}$ and define the counting function
$$\cN_{n}^{0}:=\sum_{I\in\I_{n}^{0}} \chi_{I}\:.$$
Next use John-Nirenberg inequality to get that
$$A_n^{1}:=\{x\in [0,1]\,|\, \sum_{I\in\I_{n}^{0}}\chi_{I}(x)> c\,\|\cN_n^{0}\|_{BMO_{C}}\}\:,$$
has measure less than $\frac{1}{2}$ for an appropriate choice of $c$.

Proceeding as the algorithm describes for the tiles of mass one, we repeat the above procedure and construct
\begin{itemize}
\item the collection of sets of maximal tiles $\{\p_{n}^{k,max}\}_k$;
\item the collection of sets representing the time-intervals $\{\I_{n}^{k}\}_k$;
\item the collection of counting functions $\{\cN_{n}^{k}\}_k$;
\item the level sets $\{A_n^{k}\}_k$.
\end{itemize}

Notice the following key properties of our construction:
\begin{itemize}
\item $A_n^{k}\subset A_n^{l}$ and $|A_n^{k}|\leq 2^{-(k-l)}\,|A_n^{l}|$ for any $k\geq l$;
\item $\sup_{k}\|\cN_n^{k}\|_{BMO_{C}}\leq 2^{n}$ and $\sup_{k}\|\cN_n^{k}\|_{L^{\infty}(A_n^{k}\setminus A_n^{k+1})}\lesssim 2^{n}\:;$
\item if we set $\cN_{n}:=\sum_{k} \cN_{n}^{k}$, then $\|\cN_n\|_{BMO_{C}}\lesssim 2^{n}\:.$
\end{itemize}

Now for each $k\in\N$ set
\beq\label{pk}
\p_n^k:=\left\{P=[\o,I]\,|\,\begin{array}{rl} I\subseteq A_n^{k},\:I\nsubseteq A_n^{k+1}\:\textrm{and}\:\:\:\:\:\:\:\\ A_{\P\setminus\bigcup_{j<n}\P_j,A_{n}^{k}}(P)\in [2^{-n}, 2^{-n+1})\end{array}\right\}\:,
\eeq
and let
\beq\label{Pn}
\P_n:=\bigcup_{k\geq0}\p_n^k\:.
\eeq

Conclude that
\beq\label{P}
\P=\bigcup_{n\geq0}\P_n\:.
\eeq

This ends the mass partition of our set $\P$.

In \cite{lv3}, the author shows that each family $\P_n$ may be decomposed as a union of $C\,n$ ($C\in\N$ some absolute constant) $BMO-$forests of $n^{th}$ generation. While this is not a problem if one aims for $L^p$ bounds ($1<p<\infty$), when p=1 this becomes a serious threat. Still, by carefully inspecting the proof in \cite{lv3}, one can actually be more precise: $\P_n$ may be decomposed in at most $C$ sets $\bigcup_{l=1}^{n} \P_{n,l}$ with $\{\P_{n,l}\}_l$ a ``special" sequence of $BMO-$forests that obeys: for each $l\in\N$ applying the decomposition in Definition \ref{forest} one can write $\P_{n,l}=\bigcup_{k}\P_{n,l}^{k}$  with $\P_{n,l}^{k}$ an $L^{\infty}-$forest having the $L^{\infty}$ norm of the counting function $\lesssim 2^{l}$. This fact solves the possible threat and thus, from now on, we will always suppose that $\P_n$ stands for a $BMO$ forest of $n^{th}$ generation.

\subsection{The Calderon-Zygmund decomposition - $\a$ discretization}
$\newline$

Unlike the mass decomposition, the Calderon-Zygmund decomposition presented below depends on the environment, namely on the function to which
we apply the Carleson operator $T$.

Thus, once for all, fix a nonzero function $f\in L^1$. Choose now $n\in \N$ and focus on the $BMO-$forest $\P_n$.

Our aim is to write
$$\P_n=\bigcup_{\a\in\Z}\P_n^{\a}\,,$$
with each $\P_n^{\a}$ being ``nicely" structured and, in particular, forcing \eqref{suport} to hold.


We start by decomposing  the torus in the corresponding level sets of the (dyadic) maximal function associated to $f$.
First notice that without loss of generality we may assume that $\exists\: N\in\Z$ such that
$$2^{N}<\|Mf\|_{\infty}\leq 2^{N+1}\:.$$
On the other hand since our function $f\in L^1(\TT)$ is nonzero then $$2^{M+1}\geq\frac{\int_{\TT}|f|}{|\TT|}> 2^{M}$$ for some $N\geq M\in \Z$.

Set now $\J_{\a}$ as the collection of maximal dyadic intervals $J$
such that
$$\frac{\int_{J}|f|}{|J|}> 2^{-\a}\:,$$
and let $\bar{\J}_{\a}:=\bigcup_{J\in\J_{\a}} J=\{Mf > 2^{-\a}\}\:.$

We then have that \beq\label{levset}\eeq
\begin{itemize}
\item $\forall\:\:\a\in Z,\:\:\a\leq-(N+1) \:\textrm{we have}\:\bar{\J}_{\a}=\emptyset\,;$
\item $\textrm{if}\:\:-N\leq\a_1<\a_2\leq-M\:\:\textrm{then}\:\:\bar{\J}_{\a_1}\subsetneq\bar{\J}_{\a_2}\,;$
\item $\forall\:\:\a\in Z,\:\:\a\geq-(M+1) \:\textrm{we have}\:\bar{\J}_{\a}=\bar{\J}_{-M}=\TT\,.$
\end{itemize}

For decomposing our family $\P_n$ we follow now an (increasing) inductive process.
First set
\beq\label{PN}
\P_n^{-N}:=\{P\in\P_n\,|\, \exists\: J\in \J_{-N}\: \textrm{s.t.}\: I_P\cap 51\,J\not=\emptyset\:\:\:\&\:\:\:|I_P|\leq |J|\}\,.
\eeq
Suppose that we have constructed the set $\P_{n}^{\a-1}$. Then let
$$\bar{\P}_n^{\a}:= \P_n\setminus\bigcup_{j\leq \a-1} \P_n^{j}$$ and define
\beq\label{Palfa}
\P_n^{\a}:=\left\{P\in\bar{\P}_n^{\a}\,\large|\, \exists\: J\in \J_{\a}\: \textrm{s.t.}\: I_P\cap 51\,J\not=\emptyset
\:\:\:\&\:\:\:|I_P|\leq |J|\right\}\:.
\eeq
Based on \eqref{levset}, we notice that our decomposition will end in a finite number of steps since
\beq\label{tiledec}\eeq
\begin{itemize}
\item $\forall\:\:\a\in Z,\:\:\a\leq-(N+1) \:\:\textrm{we can set}\:\:\P_n^{\a}=\emptyset\,;$
\item $\forall\:\:\a\in Z,\:\:\a\geq-(M+1) \:\:\textrm{we have}\:\:\P_n^{\a}=\emptyset\,.$
\end{itemize}
Thus, we obtained a partition of the collection $\P_n$ into
\beq\label{part}
\P_n=\bigcup_{\a\in\Z} \P_n^{\a}\;,
\eeq
such that
\beq\label{prop}
\textrm{supp}\,\P_n^{\a}:=\bigcup_{P\in \P_n^{\a}} I_P\subset 100 \,\{Mf > 2^{-\a}\}\:.
\eeq

We end this section with several important

\noindent\textbf{Observations.} 1) Remark that this partition of $P_n$ conserves the convexity property on which
the boundedness of the trees is heavily relying. More precisely we have that
\beq\label{conv}
\textrm{if}\:\:P_1<P_2<P_3\:\:\textrm{such that}\:\:P_1,\,P_3\in\P_n^{\a}\:\:\textrm{then}\:\:P_2\in\P_n^{\a}\;.
\eeq
\noindent 2) Notice that for any $P\in\P_n^{\a}\not=\emptyset$ we have that
\begin{itemize}
\item \textit{either} $\forall\:\:J\in \J_{\a-1}\:\: I_P\cap 51\,J=\emptyset\,,$
\item \textit{or} if $J\in \J_{\a-1}\: \textrm{s.t.}\: I_P\cap 51\,J\not=\emptyset\:\:\textrm{then}\:\:|I_P|> |J|$\:.
\end{itemize}
3) For $P=[\o_P,I_P]\in\P_n^{\a}\not=\emptyset$ let $c(I_P)$ be the center of the interval $I_P$ and define $I_{P^*}=[c(I_P)-\frac{17}{2}|I_P|,\,c(I_P)-\frac{3}{2}|I_P|]\cup[c(I_P)+\frac{3}{2}|I_P|,\,c(I_P)+\frac{17}{2}|I_P|]$. Then we have that
\beq\label{support}
\textrm{supp}\,T_{P}\subseteq I_P\:\:\:\textrm{and}\:\:\:\textrm{supp}\,T_{P}^{*}\subseteq I_{P^*}\:.
\eeq
Moreover writing $$I_{P^*}=\bigcup_{r=1}^{14}I_{P*}^r$$ with each $I_{P^*}^r$ a dyadic
interval of length $|I_P|$ we have the following property
\beq\label{suppPstar}
\frac{\int_{I_{P^*}^r}|f|}{|I_P|}< 2^{-\a+10}\:.
\eeq
4) Let $\p\in\P_n^{\a}$ be a tree. Define $\p^{min}$ the collection of minimal tiles in $\p$. Further set
$$\J_{\p^*}:=\{I_{P*}^r\,|\,P\in\p^{min}\:\:\textrm{and}\:\:r\in\{1,\ldots,14\}\}\;,$$
and $CZ(\J_{\p^*})$ the Calderon-Zygmund decomposition of the interval $[0,1]$ with respect to $\J_{\p^*}$.

Then, from \eqref{Palfa} and \eqref{suppPstar} we deduce the following \textbf{key property}:
\beq\label{cz}
\frac{\int_{I}|f|}{|I|}< 2^{-\a+10}\:\:\:\:\:\:\:\:\forall\:I\in CZ(\J_{\p^*})\:.
\eeq

\section{The proof of the Main Theorem.}

\subsection{Proof of a)} Our aim is to show that with the notions previously defined we have
\beq\label{a}
\textrm{supp}\,T^{\P_n^{\a}}\subseteq 100\,\{Mf>2^{-\a}\}\:\:\:\textrm{and}\:\:\:\|T^{\P_n^{\a}}\,f\|_1\lesssim 2^{-\a}\,|\{Mf>2^{-\a}\}|\:.
\eeq
The first of the above conditions is an immediate consequence of the construction of the tile families $\{\P_n^{\a}\}$\:.

For the second condition we need to analyze the structure of each $\P_n^{\a}$.
Let $\P_n^{\a}=\bigcup_{j\in\N} \P_n^{\a,j}$ be the decomposition of $\P_n^{\a}$ in $L^{\infty}-$forests. Further, for $j\in\N$, we decompose
each $L^{\infty}-$forest $\P_n^{\a,j}=\bigcup_{k}\p_{k}^{\a,n,j}$ in maximal trees. Set $I_{P_{k}^{\a,n,j}}$ the time interval of the top of the tree $\p_{k}^{\a,n,j}$ and define $I_{\P_n^{\a,j}}:=\bigcup_{k}I_{P_{k}^{\a,n,j}}$ and $I_{\P_n^\a}:=\bigcup_{j}I_{\P_{n}^{\a,j}}$ respectively.

Now, as a consequence of the construction in Section 5, for $\cN_j^{\a,n}:=\sum_{k}\chi_{I_{P_{k}^{\a,n,j}}}$ and $\cN^{\a,n}=\sum_{j}\cN_j^{\a,n}$,
we have
\beq\label{counting}
\|\cN_j^{\a,n}\|_{\infty}\lesssim 2^{n}\:\:\:\:\:\:\:\forall\:j\in\N\:\:\textrm{and}\:\:\:\:\:\|\cN^{\a,n}\|_1\lesssim 2^{n}\,|I_{\P_n^{\a}}|\:.
\eeq
Based on \eqref{counting}, it is thus enough to prove that for $\p\subset\P_n^{\a}$ tree
\beq\label{tree}
\int|T^{\p}f|\lesssim 2^{-n}\,2^{-\a}\,|I_{\p}|\;.
\eeq

Indeed, from \eqref{tree} and \eqref{counting} we deduce
 \beq\label{fin}
\begin{array}{cl}
\int|T^{\P_n^\a}f|\lesssim 2^{-n}\,2^{-\a}\,\|\cN^{\a,n}\|_1\,\lesssim 2^{-\a}\, |I_{\P_n^{\a,0}}|\lesssim 2^{-\a}\,|\{Mf>2^{-\a}\}|\;.
\end{array}
\eeq

Now for showing \eqref{tree} we proceed as follows:

Without loss of generality we may assume that all the tiles $P\in\p$ are at a constant frequency $\o$ (the frequency of the tree).

Define
\beq\label{proj}
\LL_{\p}(f):=\sum_{J\in CZ(\J_{\p^*})} \frac{\int_{J}f(s)\,e^{i\,\o\,s}\,ds}{|J|}\chi_J\;.
\eeq

Observe that as a consequence of \eqref{cz} we have
\beq\label{linfproj}
|\LL_{\p}(f)|\lesssim 2^{-\a}\,\chi_{I_{\p}}\:.
\eeq

Setting now $\p_0$ the shift of $\p$ to the origin, we have
$$\int|T^{\p}f|\leq \int|T^{\p_0}(e^{i\,\o\,\cdot}f(\cdot)-\LL_{\p}(f)(\cdot))|\,+\,\int|T^{\p_0}(\LL_{\p}(f)(\cdot))|\:. $$

For the first term we use the mean zero condition
$$\int_{I_{P^*}}\{f(\cdot)\,e^{i\,\o\,\cdot}-\LL_{\p}(f)(\cdot)\}=0\:\:\:\:\:\:\:\:\forall\: P\in\p$$
 and deduce that
\beq\label{t1}
\int|T^{\p_0}(e^{i\,\o\,\cdot}f(\cdot)-\LL_{\p}(f)(\cdot))|\lesssim 2^{-n}\,\|f\|_{L^1(I_{\p})}\:.
\eeq

For the second term, we use relation \eqref{linfproj} and the $L^2$ boundedness of the Hilbert transform
\beq\label{t2}
\int|T^{\p_0}(\LL_{\p}(f))|\lesssim |E(\p)|^{\frac{1}{2}}\,2^{-n/2}\,\|\LL_{\p}(f)\|_{L^2(I_{\p})}\lesssim 2^{-n}\,2^{-\a}\,|I_{\p}|\;,
\eeq
where we set $E(\p):=\bigcup_{P\in\p} E(P)$.

Here we have used the key Carleson measure estimate
\beq\label{Carlmeas}
|E(\p)|\lesssim 2^{-n}\,|I_{\p}|\,.
\eeq

Indeed, for proving \eqref{Carlmeas}, we follow the reasoning from \cite{lv1} (p. 481) and for
$$\J_{\p}:=\{I_{P}\,|\,P\in\p^{min}\}\,,$$ we set
\begin{align*}
\check{\J}(\p):=\big\{I\subset I_{\p}\:|\: & \textrm{Exactly one of the left or right halves}\\
&\textrm{of $I$ contains an element of $\J_{\p}$} \big\} \: \cup \J_{\p}\:.
\end{align*}
and $\breve{\p}=\left\{P=[\o,\,I]\in \p\:\:|\:\:I\in\check{\J}(\p)\:\right\}.$

Then we have
$$|E(\p)|\leq\sum_{P\in \breve{\p}}|E(P)|\lesssim 2^{-n}\:\sum_{I\in \check{\J}(\p)}|I|\lesssim 2^{-n}\:|I_{\p}|\:.$$

This ends our proof.

Notice that \eqref{a} implies
$$\|T^{\P_n}\|_1\lesssim\sum_{\a}\|T^{\P_n^{\a}}\,f\|_1\lesssim\|Mf\|_1\lesssim\|f\|_{L\log L}\:.$$

\subsection{Proof of b)} In this section we will show a slightly stronger statement than the one claimed in our theorem. More precisely,
we prove that the operators $\{T^{\P_n}\}_n$ are \emph{uniformly} weak $(1,1)$ bounded:

We claim that $\forall\:G\subseteq [0,1]\:\:\exists\:\:G'\subseteq G,\:\:|G'|>\frac{1}{2}|G|$ such that\footnote{A similar statement for a rougher
mass-discretized family $\P_n$ is proved in \cite{LaDo}. For this, the authors are using the mass-size decomposition technique presented in \cite{lt3}. Then, they interpolate the resulting $L^1$-estimate with a ``modified" $L^2$ estimate to get good control near $L^1$ - this reasoning is also used in our approach.}
\beq\label{lweak1}
\int_{G'}|T^{\P_n}\,f|\lesssim\|f\|_1\:\:\:\forall\:n\in\N\,.
\eeq

Fix $G\subset\TT$ and define $G'$ as
$$G'=\{x\in G\,|\,Mf(x)\leq C\,\frac{\|f\|_1}{|G|}\}\:,$$
where here $C>0$ is some large constant appropriately chosen. This assures the requirement $G'\subset G$ with $|G'|\gtrsim |G|$.

Set now $\ll=C\frac{\|f\|_1}{|G|}.$ It will be enough to prove that\footnote{For two positive quantities $A,\,B$ we write $A\approx B$ iff $2^{-10}\,B<A<2^{10}\,B$.}
\begin{l1}\label{alfak}
Let $k\in \N$ and suppose $2^{-\a}\approx \ll\,2^{-k}$.
Then the following relation holds:
\beq\label{kdecay}
\int_{G'}|T^{\P_n^\a}\,f|\lesssim 2^{-\frac{k}{2}}\,\|f\|_1\:\:\:\forall\:n\in\N\,.
\eeq
\end{l1}

This result relies on the tree estimate provided below

\begin{l1}\label{treek}
Let $\p\subset\P_n^\a$ be a tree with top $I_{\p}$. Then we have
\beq\label{essential}
\int_{G'}|T^{\p}\,f|\lesssim 2^{-\a}\,|G'\cap E(\p)|^{\frac{1}{2}}\,2^{-n/2}\,|I_{\p}|^{\frac{1}{2}}.
\eeq
\end{l1}

\begin{proof}

Let $CZ(\J_{\p^*})$ be the Calderon-Zygmund decomposition described in the Observations, Section 5.2..  As before, we can assume\footnote{To set our problem in the context of this assumption we take advantage of the translation invariance of our statement and use a standard estimate (see \cite{f}) that gives an error of order $2^{-\a}\,|G'\cap E(\p)|$.}
without loss of generality that all the tiles $P\in\p$ are at a constant frequency $\o=0$.

Now following the description from a) we have
$$\int_{G'}|T^{\p}f|\leq \int_{G'}|T^{\p}(f-\LL_{\p}(f))|\,+\,\int_{G'}|T^{\p}(\LL_{\p}(f))|\:. $$

The second term is trivially bounded by the Cauchy-Schwarz inequality and the $L^2$ boundedness of the Hilbert transform
\beq\label{hilbl2}
\begin{array}{cl}
\int_{G'}|T^{\p}(\LL_{\p}(f))|\lesssim |G'\cap E(\p)|^{\frac{1}{2}}\,2^{-n/2}\,\|\LL_{\p}(f)\|_{L^2(I_{\p})}\\
\lesssim 2^{-\a}\,|G'\cap E(\p)|^{\frac{1}{2}}\,2^{-n/2}\,|I_{\p}|^{\frac{1}{2}}.
\end{array}
\eeq

For the first term we need to be more careful; we will show that for any $g\in L^{\infty}(\TT)$ with $\textrm{supp}\:g\subseteq G'$ we have
\beq\label{hilblinfty}
\left|\int({T^{\p}}^{*}g)\,(f-\LL_{\p}(f))\right|\lesssim 2^{-\a}\,|G'\cap E(\p)|^{\frac{1}{2}}\,2^{-n/2}\,|I_{\p}|^{\frac{1}{2}}\,\|g\|_{\infty}.
\eeq
At this point we make essential use of the mean zero property $$\int_{J}(f-\LL_{\p}(f))=0\:\:\:\:\forall\: J\in CZ(\J_{\p^*})\:.$$

Thus, for proving \eqref{hilblinfty}, it is enough to show that for $g=g\chi_{G'}\in L^{\infty}$
\beq\label{hilblinfty1}
\int\left|({T^{\p}}^{*}g-\LL_{\p}({T^{\p}}^{*}g))\,(f-\LL_{\p}(f))\right|\lesssim
2^{-\a}\,|G'\cap E(\p)|^{\frac{1}{2}}\,2^{-n/2}\,|I_{\p}|^{\frac{1}{2}}\,\|g\|_{\infty}.
\eeq

For fixed $J\in CZ(\J_{\p^*})$  and $x\in J$ we have
$$\left|{T^{\p}}^{*}g(x)-\frac{1}{|J|}\int_{J}{T^{\p}}^{*}g(s)ds\right|=$$
\beq \label{est}
\left|\frac{1}{|J|}\int_{J}\left\{\sum_{{P\in\p}\atop{|I_P|\geq
|J|}}\int_{\TT}\left[\f_k(x-y)-\f_k(s-y)\right]g(y)\chi_{E(P)}(y)dy\right\}ds\right|\lesssim
\eeq
$$\frac{1}{|J|}\int_{J}\left\{\sum_{I_{P^*}\supseteq J} |I_P|^{-1}\,|J|\frac{\int_{E(P)\cap G'}|g|}{|I_P|}\right\}ds\lesssim
\|g\|_{\infty}\,\frac{\chi_{J}}{|J|}\sum_{I_{P^*}\supseteq J}\frac{|J|^2}{|I_P|^2}\, |G'\cap E(P)|\;.$$

This last relation gives us
$$\int \left|({T^{\p}}^{*}g-\LL_{\p}({T^{\p}}^{*}g))\,(f-\LL_{\p}(f))\right|$$
$$\lesssim 2^{-\a}\,2^{-\frac{n}{2}}\,\|g\|_{\infty}\,\sum_{J\in CZ(\J_{\p^*})}
\sum_{I_{P^*}\supseteq J}|J|^{\frac{1}{2}}\,(\frac{|J|}{|I_P|})^{\frac{3}{2}}\,|G'\cap E(P)|^{\frac{1}{2}}$$
$$\lesssim 2^{-\a}\,|G'\cap E(\p)|^{\frac{1}{2}}\,2^{-n/2}\,|I_{\p}|^{\frac{1}{2}}\,\|g\|_{\infty}\:,$$
thus proving \eqref{essential}.

Here we have relayed on the following key relation:
\beq \label{CMG}
\sum_{{P\in\p}}\sum_{{J\in CZ(\J_{\p^*})}\atop{J\subseteq I_{P^*}}}(\frac{|J|}{|I_P|})^{\frac{3}{2}}\,|J|^{\frac{1}{2}}\,|G'\cap E(P)|^{\frac{1}{2}}
\lesssim |G'\cap E(\p)|^{\frac{1}{2}}\,|I_{\p}|^{\frac{1}{2}}\,.
\eeq

For proving this, take $l\in\N$ and set
$$CZ_{l}^{*}(I_P):=\{J\in CZ(\J_{\p^*})\,|\,J\subseteq I_{P^*}\,,\:\:\:|J|\approx 2^{-l}\,|I_P|\,\}\:.$$

Now \eqref{CMG} will be a consequence of
\beq \label{lcz}
S_l:=\sum_{{P\in\p}}\sum_{J\in CZ_{l}^{*}(I_P)} |J|^{\frac{1}{2}}\,|G'\cap E(P)|^{\frac{1}{2}}
\lesssim_{\ep} 2^{l(\frac{1}{2}+\ep)}\,|G'\cap E(\p)|^{\frac{1}{2}}\,|I_{\p}|^{\frac{1}{2}}\:,
\eeq
where here $\ep\in (0,1)$ is some absolute constant.

For proving \eqref{lcz}, we start by refining the set $\p$ as follows:

For $u\in\{0\,\ldots,l\}$ we set $$\p^{u}:=\{P\in\p\,|\,2^{u-1}<\#CZ_{l}^{*}(I_P)\leq 2^{u}\}\:.$$

Then, fix $u$ and define
$$S_{l,u}:=\sum_{{P\in\p^{u}}}\sum_{J\in CZ_{l}^{*}(I_P)} |J|^{\frac{1}{2}}\,|G'\cap E(P)|^{\frac{1}{2}}\:.$$

Next, we decompose inductively the set $\p^{u}$ as follows: set $\p_1$ the collection of maximal (with respect to $``\leq"$) tiles in $\p^{u}$; repeat this procedure for the collection $\p\setminus \p_1$ and thus construct the set of maximal tiles $\p_2$. Continue this process until
exhausting $\p^{u}$. Thus we end up with partitioning
$$\p^u=\bigcup_{r=1}^m \p_r\,,$$
into collections of (successively) maximal tiles.

Applying now Cauchy-Schwarz we deduce that
\beq \label{cs}
S_{l,u}\leq |G'\cap E(\p)|^{\frac{1}{2}}\,\sum_{r=1}^{m} \left(\sum_{P\in\p_r}(\sum_{J\in CZ_{l}^{*}(I_P)}|J|^{\frac{1}{2}})^2\right)^{\frac{1}{2}}\:.
\eeq

Finally, we observe that
\beq \label{carlmeas}
\sum_{r=1}^{m} \left(\sum_{P\in\p_r}(\sum_{J\in CZ_{l}^{*}(I_P)}|J|^{\frac{1}{2}})^2\right)^{\frac{1}{2}}\lesssim 2^{\frac{l}{2}}\,|I_{\p}|^{\frac{1}{2}}\:.
\eeq

This last fact is a consequence of the construction of the sets $\{\p_r\}_r$, the definition of $CZ_{l}^{*}(I_P)$ and of the fact that for any
 $J\in CZ(\J_{\p^*})$ one has
$$\#\{P\in\p\,|\,J\in CZ_{l}^{*}(I_P)\}\lesssim 1\:.$$

Indeed, based on these above mentioned facts, one has that the main contribution in the left hand side of \eqref{carlmeas} comes from the first $O(2^{l-u})$ terms of the sum. We leave further details for the interested reader.

\end{proof}

We pass now to the proof of Lemma \ref{alfak}. With the notations from point a) and based on \eqref{essential}, for each $L^{\infty}$-forest $\P_n^{\a,j}$  we decompose $\P_n^{\a,j}=\bigcup_{k}\p_{k}^{\a,n,j}$ into maximal trees and deduce
$$\int_{G'}|T^{\P_{n}^{\a,j}}\,f|\lesssim\sum_{k\in\N}\int_{G'}|T^{\p_{k}^{\a,n,j}}\,f|\lesssim 2^{-\a}\,\sum_{k\in\N}|G'\cap E(\p_{k}^{\a,n,j})|^{\frac{1}{2}}\,2^{-n/2}\,|I_{\p_{k}^{\a,n,j}}|^{\frac{1}{2}}$$
$$\leq2^{-\a}\, 2^{-n/2}\,\{\sum_{k\in\N}|G'\cap E(\p_{k}^{\a,n,j})|\}^{\frac{1}{2}}\,\{\sum_{k\in\N}|I_{\p_{k}^{\a,n,j}}|\}^{\frac{1}{2}}$$
$$\lesssim 2^{-\a}\,|G'\cap E(\P_{n}^{\a,j})|^{\frac{1}{2}}\,|I_{\P_n^{\a,j}}|^{\frac{1}{2}}\:.$$
From the Carleson measure condition imposed in Definition \ref{forest} ii), we have
$$\int_{G'}|T^{\P_{n}^{\a}}\,f|\lesssim\sum_{j}\int_{G'}|T^{\P_{n}^{\a,j}}\,f|
\lesssim 2^{-\a}\,|G'|^{\frac{1}{2}}\,\sum_{j}|I_{\P_n^{\a,j}}|^{\frac{1}{2}}\lesssim
2^{-\a}\,|G'|^{\frac{1}{2}}\,|I_{\P_n^{\a,0}}|^{\frac{1}{2}}\:.$$
Finally, using that
\beq\label{max}
|I_{\P_n^{\a,0}}|\lesssim |\bar{\J}_{\a}| \approx 2^{\a} \int_{\bar{\J}_{\a}}|f|\:,
\eeq
and setting $2^{-\a}\approx \ll\,2^{-k}$ we conclude that \eqref{kdecay} holds.

\subsection{Proof of c)} The central estimate for proving \eqref{lp} is given by
\beq\label{ll2}
\|T^{\P_n}\,f\|_{2}\lesssim n^2\,2^{-n/2}\,\|f\|_{2}\:.
\eeq

We limit ourselves to just providing the main ideas. For the details of the proof the reader should consult \cite{lv3}.

One first decomposes $\P_n=\bigcup \P_n^j$ with each $\P_n^{j}$ an $L^{\infty}-$forest, and then proves that
\begin{itemize}
\item $\|T^{\P_n^j}\,f\|_{2}\lesssim n\,2^{-n/2}\,\|f\|_{2}\:.$
\item for each $k\in\{0,\ldots,n-1\}$ the sequence $\{T^{\P_n^{k+j\,n}}\}_{j\in\N}$ consists from almost orthogonal operators
and thus one can apply Cotlar-Stein lemma.
\end{itemize}

This ensures that \eqref{ll2} holds. Now, for $1<p\leq 2$, relation \eqref{lp} is just a consequence of \eqref{l1}, \eqref{ll2} and classical interpolation theory. For $p>2$ we can proceed as in the $L^2$ case: one first proves the desired estimate for an $L^\infty$ forest and then using the structure of the $BMO$-forest one is able to extend the initial result to the entire family $\P_n$ (for more details see \cite{lv3}).

\subsection{Proof of d)} We will show\footnote{In \cite{ArB}, one can found (without a proof) the following statement: For $1\leq p\leq 2$ we have that $\|Tf\|_{p,\infty}\lesssim \|f\|_p\,\log\frac{e\,\|f\|_{2}}{\|f\|_p}\:.$ The unpublished proof (\cite{Arunp}) of this result handled to me by Arias de Reyna, relies on Carleson's original approach to the pointwise convergence of the Fourier Series. Finally, one may notice that the case $p=1$ in his result is equivalent with the case $p=2$ in \eqref{12}.} that for $1<p<\infty$
\beq\label{12}
\|Tf\|_{1,\infty}\lesssim_p \|f\|_1\,\log\frac{e\,\|f\|_{p}}{\|f\|_1}\:.
\eeq

As before, we start by reformulating \eqref{12} in the equivalent form

$\forall\:G\subset\TT\:\:\exists\:\:G'\subset G$ with $|G'|\gtrsim |G|$ such that
\beq\label{ref12}
\int_{G'} |Tf|\lesssim_p \|f\|_1\,\log\frac{e\,\|f\|_{p}}{\|f\|_1}\:.
\eeq

Repeat the same construction as in a). Then, with the same notations from c), we have
\beq\label{forestl2}
\int_{G'}|T^{\P_n}f|\lesssim _p 2^{-\d\,n/p^*}\,\|f\|_{p}\:.
\eeq
Combining \eqref{l1} and \eqref{forestl2} we conclude that
$$\int_{G'} |Tf|\leq \sum_{n\in\N}\int_{G'} |T^{\P_n}f|\lesssim_p \sum_{n\in\N} \min\{ \|f\|_1,\,2^{-\d\,n/p^*}\,\|f\|_p\}\lesssim_p\|f\|_1\,\log\frac{e\,\|f\|_{p}}{\|f\|_1}\,,$$
proving the desired result.

\section{Remarks}

1) As mentioned in the introduction, the central part of our paper resides on the tile decomposition described in Section 5.

The mass discretization of $Tf$ is independent\footnote{The mass parameter depends on the function $N$ which may be taken as just an arbitrary measurable function as long as the final estimates on the operator $T$ do not depend on it.} of $f$ and results in the geometric organization of $\P=\bigcup_n \P_n$ which
gives us both the $n$-decay\footnote{This helps us in both summing the operators $\{T^{\P_n}\}_n$ (and thus obtaining $L^p$ bounds for $1<p<\infty$) and in interpolating with norm estimates near $L^1$.} for each $\{\|T^{\P_n}f\|_p\}_n$ with $1<p<\infty$ ($f\in L^p$) and the ability to sum\footnote{Here is the point where our technique overcomes the difficulty of treating the exceptional sets.} within each ``scale" $n$ the lengths of the time support of the maximal trees in $\P_n$.

The Calderon-Zygmund discretization $\P_n=\bigcup_{\a}\P_n^{\a}$ realizes the decomposition of the function $T^{\P_n}f$ depending on $f$ and is designed to get a good control near $L^1$. This discretization accounts for the multi-frequency nature of our problem. A similar instance appeared in \cite{NOT}, where the authors are elaborating a ``Calderon-Zygmund decomposition for multiple frequencies" for a given function $f$. Our approach though is quite different: instead of decomposing the input object (the function $f$) at multiple frequencies imposing a mean zero condition of the initial function for each frequency, we rely on the properties\footnote{Here the algorithm described for the mass decomposition plays a fundamental role.} of the initial discretization $\P=\bigcup_n \P_n$ and first decompose the tile family $\P_n$ in subfamilies $\{\P_n^{\a}\}_{\a}$ followed by a further decomposition of the corresponding output object $T^{\P_n^{\a}}f$ in multiple pieces with each piece having the mean zero condition strictly relative to the frequency at which it lives.

In the setting of the present paper our procedure is more effective and unlike the one in \cite{NOT} gives an explicit construction. On the other hand, the decomposition in \cite{NOT} is more general and hence can be used in other problems which are not necessarily related to the Carleson operator.
$\newline$

2) The treatment of the ``forest" operators $T^{\P_n}$ in our Theorem offers a substitute to the classical theory of the Calderon-Zygmund operators:
indeed, following the tile decomposition and the proofs of a) and b) one notices that the weak $(1,1)$ bound is obtained by using a Vitali type covering argument - in which we can sum the lengths of the intervals of the trees in the structure of $T^{\P_n}f$ depending on the size of the maximal function associated to $f$. Also like in the classical theory, through our decomposition the entire weight moves on the $L^2$ behavior of $T^{\P_n}$ where one needs to use orthogonality methods. To complete our parallelism, it would be of interest if our methods could provide a satisfactory theory for the adjoint operator ${T^{\P_n}}^{*}$ near $L^1$.

We think this topic deserves further investigation.
$\newline$

3) Finally, in view of our approach, the following question appears as natural:
$\newline$

\noindent\textbf{Open question.} \textit{Fix $f\in L^1(\TT)$. With the previous notations and definitions, set
$$\P^{\a}:=\bigcup_{n}\P_n^{\a}\,,\:\:\:\:\textrm{with}\:\:\:\a\in\N\,.$$
\indent Is it then true that $\exists\:C>0$ absolute constant such that for any $\a\in\N$
\beq\label{oq}
\|T^{\P^{\a}}\,f\|_{1}\leq C \|f\|_1\,?
\eeq}

A positive answer to this question would imply that
\beq\label{Carl}
\|Tf\|_{L^1}\lesssim \|f\|_{L\log L}\,,
\eeq
which is the best one can hope for if we require strong $L^1$ bounds for the Carleson operator $T$.

Of course it would be still very interesting if \eqref{oq} holds with the $L^1$ norm replaced by the $L^{1,\infty}$ norm in the left hand side.

With the current technology we can only prove that
 $$\|T^{\P^{\a}}\,f\|_{1}\leq C \|f\|_{L \log L}\,.$$

\section{Appendix - spaces near $L^1$}

In this section we briefly introduce the definitions of the relevant rearrangement invariant Banach spaces which appeared in the previous literature
when studying the problem of the pointwise convergence of the Fourier Series near $L^1$.
\begin{d0}\label{orlicz}
Let $\f:\:[0,\infty)\,\mapsto\, [0,\infty)$ be an absolutely continuous function with the following properties:
\begin{itemize}
\item $\exists\,C>0$ such that $\f(t^2)\leq C\,\f(t)$ for all $t\geq 0$;
\item $\f'(t)\geq 0$ almost everywhere;
\item $\f(0)=0$;
\item $\lim_{t\rightarrow \infty} \f(t)=\infty.$
\end{itemize}
Then we define the space $L\,\f(L)$ of all (measurable) functions $f$ for which
\beq\label{O}
\|f\|_{\l\,\f(\l)}:=\int_\TT |f|\,|\f(f)|<\infty\:.
\eeq
\end{d0}

A classical result in Banach space theory (see e.g. \cite{BS}) asserts

\noindent\textbf{Proposition.} \textit{The space $L\,\f(L)$ endowed with the norm
\footnote{Here $f^*$ stands for the decreasing rearrangement of $f$.}
\beq\label{norm}
\|f\|_{L\,\f(L)}:=\int_\TT f^*(t)\,\f(f^{*})(t)\,dt<\infty\:,
\eeq
becomes a Banach rearrangement space.}

In the topic treated in this paper three special choices for the function $\f$ are of interest:

\begin{itemize}
\item $\f(t)= \log(1+t)$ - defines the Zygmund space $L\,\log L$;

\item $\f(t)= \log(1+t)\,\log\log(1+t)$ - defines the space $L\,\log L\,\log\log L$ considered by Sj\"olin in \cite{sj3}.

\item $\f(t)= \log(1+t)\,\log\log\log (10+t)$ - defines the space $L\,\log L\,\log\log\log L$ considered by Antonov in \cite{An}.
\end{itemize}

To complete the picture and thus present the evolution until nowadays, we need to consider two more spaces:

The first space was considered by F. Soria in \cite{So1}:

\noindent Let $B_{\f}$ be the set of the measurable functions for which
$$\|f\|_{\f}:=\int_{0}^{\infty} \f (\lambda_f(t))\,dt<\infty$$ where
here $\lambda_f$ is the distribution function of $f$ given by $\lambda_f(t)=|\{x\in\TT\,|\,|f(x)|>t\}|$. Take now the subspace
$B_{\f}^*\subset B_{\f}$ defined by $$B_{\f}^*:=\{f\,|\,\|f\|_{\f}^*<\infty\}$$ where
$\|f\|_{\f}^{*}=\int_{0}^{\infty}\f(\lambda_f(t))(1+\log(\frac{\|f\|_{\f}}{\f(\lambda_f(t))}))dt.$

The pointwise convergence theory developed by Soria was addressing the space $B_{\f_1}^*$ with $\f_1(s)= s(1+\log^{+}\frac{1}{s})$.

The second space was introduced by Arias de Reyna in \cite{Ar}:

\noindent Let $QA$ be the quasi-Banach space defined as follows:
$$QA:=\{f:\:\TT\mapsto C\,|\,f\:\textrm{measurable},\:\|f\|_{QA}<\infty \}\:\:\:\textrm{where}$$
$$\|f\|_{QA}:=\inf\left\{\sum_{j=1}^{\infty}(1+\log j)\|f_j\|_1\,\log\frac{e\,\|f_j\|_{\infty}}{\|f_j\|_1}\:\:\left|\right.\:\:
\begin{array}{cl}
f=\sum_{j=1}^{\infty}f_j,\:\\
\sum_{j=1}^{\infty}|f_j|<\infty\:\textrm{a.e.}
\end{array}
\right\}\;.$$

 With the exception of Zygmund's $L\log L$ space, it is known (\cite{sj3}, \cite{An}, \cite{So1} and \cite{Ar}) that all the other spaces considered here, namely $L\,\log L\,\log\log L$, $L\,\log L\,\log\log\log L$, $B_{\f_1}^*$ and $QA$  are rearrangement invariant spaces of functions with almost everywhere convergent Fourier series.

\end{document}